\documentclass[12pt]{article}
\usepackage{amsmath, amsfonts, amssymb, amscd, enumerate}
\usepackage[dvips]{color}
\usepackage{graphicx}

\def\a{\alpha}
\def\b{\beta}

\def\i{\iota}

\def\o{\omega}

\def\u{\upsilon}

\chardef\tempcat=\the\catcode`\@
\catcode`\@=11
\def\cyracc{\def\u##1{\if \i##1\accent"24 i
    \else \accent"24 ##1\fi }}
\newfam\cyrfam

\DeclareFontFamily{OT1}{msb}{}{}
\DeclareFontShape{OT1}{msb}{m}{n}
 {  <5> <6> <7> <8> <9> <10> gen * msbm
      <10.95><12><14.4><17.28><20.74><24.88>msbm10}{}
\DeclareMathAlphabet{\bubble}{OT1}{msb}{m}{n}

\def\bR{{\mathbb R}}

\def\bC{{\mathbb C}}
\def\bH{{\mathbb H}}

\def\bO{{\mathbb O}}

\newfont{\goth}{eufm10 scaled \magstep1}
\def\ga{\mbox{\goth a}}

\def\gf{\mbox{\goth f}}
\def\gg{\mbox{\goth g}}
\def\ge{\mbox{\goth e}}
\def\gh{{\mbox{\goth h}}}
\def\gk{\mbox{\goth k}}
\def\gl{\mbox{\goth l}}
\def\gm{\mbox{\goth m}}
\def\gn{\mbox{\goth n}}

\def\gp{\mbox{\goth p}}
\def\gq{\mbox{\goth q}}

\def\gu{\mbox{\goth u}}

\def\gz{\mbox{\goth z}}

\def\gso{\mbox{\goth so}}
\def\gsu{\mbox{\goth su}}
\def\gspin{\mbox{\goth spin}}
\def\gsp{\mbox{\goth sp}}
\def\gsl{\mbox{\goth sl}}
\def\ggl{\mbox{\goth gl}}

\newfont{\mcal}{eusm10 scaled \magstep1}

\def\exp{\mathrm{exp\;}}
\def\ra{\rightarrow}
\def\Id{\mathrm{Id\;}}
\def\tr{\mathrm{tr\;}}

\def\Ker{\mathrm{Ker\;}}

\def\Ad{\mathrm{Ad}}
\def\ad{\mathrm{ad}}

\def\wt{\widetilde}

\def\exp{ \mathrm exp\;}
\def\span{ \mathrm span \;}
\newtheorem{Th}{Theorem}
\newtheorem{Prop}{Proposition}
\newtheorem{Cor}{Corollary}
\newtheorem{Lem}{Lemma}
\newtheorem{Def}{Definition}
\newtheorem{Rem}{Remark}
\def\bt{\begin{Th}}
\def\et{\end{Th}}
\def\bp{\begin{Prop}}
\def\ep{\end{Prop}}
\def\bc{\begin{Cor}}
\def\ec{\end{Cor}}
\def\bl{\begin{Lem}}
\def\el{\end{Lem}}
\def\bd{\begin{Def}}
\def\ed{\end{Def}}
\def\pf{\noindent{\it Proof. }}
\def\qed{\hspace{2ex} \hfill $\square $ \par \medskip}

\def\be{\begin{equation}}
\def\ee{\end{equation}}
\def\arr{\begin{array}{rlll}}
\def\ea{\end{array}}
\def\bea{\begin{eqnarray}}
\def\eea{\end{eqnarray}}
\def\bean{\begin{eqnarray*}}
\def\eean{\end{eqnarray*}}

\def\span {\mathrm {span} \, }
\def\nbh {neighborhood \, }
\def\hpsRms{ homogeneous pseudo-Riemannian manifolds \,}
\def\hpsRm{ homogeneous pseudo-Riemannian manifold\,}
\def\hRm{ homogeneous Riemannian manifold\,}
\def\dRhm {degenerate Riemannian homogeneous  manifold\,}
\def\psR{pseudo-Riemannian \,}

\def\hLm { homogeneous Lorentzian  manifold \,}
\def\hLms { homogeneous Lorentzian  manifolds \,}
\def\Lm {Lorentzian manifold}
\def\C1Lm { cohomogeneity one Lorentzian manifold \,}

\def\Lie{ \mathrm {Lie} \,}
\def\diag{ \mathrm {diag} \,}
\definecolor{red}{rgb}{1,0,0}           
  \definecolor{green}{rgb}{0,1,0}
  \definecolor{blue}{rgb}{0,0,1}

\begin{document}

\begin{titlepage}
 \rightline{}
\vskip 1.5 true cm
\begin{center}
{\Large Homogeneous Lorentzian manifolds of semisimple group}
\end{center}
\vskip 1.0 true cm
\begin{center}
{  D.V.\ Alekseevsky } \\[3pt]
\end{center}
\vskip 2.0 true cm \baselineskip=18pt

\centerline{Edinburgh University   and  Hamburg University}
\vskip 0.5 truecm
\abstract  We   describe the  structure of  $d$-dimensional
homogeneous Lorentzian $G$-manifolds $M=G/H$ of a  semisimple  Lie
group $G$. Due to a result by N. Kowalsky, it is sufficient  to
consider  the case when the group  $G$ acts properly, that is the
stabilizer $H$ is compact.  Then  any  homogeneous space $G/\bar H$
with a smaller group $\bar H \subset H$  admits an invariant
Lorentzian metric. A homogeneous  manifold  $G/H$ with a connected
compact stabilizer $H$ is called a minimal admissible manifold  if
it admits an invariant Lorentzian metric, but no homogeneous
$G$-manifold $G/\tilde H$ with a larger connected compact stabilizer
$\tilde H \supset H$ admits such a metric.
 We give a description of  minimal homogeneous Lorentzian
 $n$-dimensional
 $G$-manifolds $M = G/H$ of a simple (compact or noncompact) Lie group $G$.
For $n \leq 11$,  we  obtain a list of all such manifolds $M$ and
describe invariant Lorentzian metrics on $M$.


\end{titlepage}
\tableofcontents
\vskip 0.5 truecm
\newpage

\section{Introduction}

  We   discuss the problem of classification of homogeneous
  Lorentzian $G$-manifolds $M = G/H$ of a semisimple Lie group $G$.
   We say that  a $G$-manifold $M$ is proper if the  action of  the
   isometry group $G$ on $M$ is proper.  In contrast with  the Riemannian
   case,  there  are nonproper homogeneous Lorentzian  manifolds,
   for example, the De Sitter  space $dS^n = SO_{1,n}/SO_{1,n-1}$ and
    the anti De Sitter  space $AdS^n = SO_{2,n}/SO_{1,n-1}$.\\
     A surprising result by
    Nadin  Kowalsky shows  that these  spaces of constant curvature
    exhaust    all nonproper homogeneous Lorentzian manifolds of a
    simple group $G$ (up to a  local isometry).\\
     This  result  had
    been generalized by M.\ Deffaf, K.\ Melnick  and A.\ Zeghib to
    the case of a
     semisimple group $G$ :\\
     Any nonproper  homogeneous Lorentzian manifold of a semisimple
     Lie group $G$ is a local product of the  (anti) De Sitter space
     and a Riemannian homogeneous manifold.\\
     This reduces  the  classification of  homogeneous
     Lorentzian  manifolds $M=G/H$ of a semisimple Lie group  to the
     case  when the stabilizer  $H$ is compact.\\
      We  will always assume  that all considered  Lie groups  are
      connected. In particular, by a stability subgroup of  an
      action of a Lie group on a manifold  we will understand  a
      connected  stability subgroup.

      We  say  that a proper homogeneous manifold $M=G/H$ ( and
       the  stability  subgroup  $H$ ) is
      { \bf admissible } if $M$ admits an invariant Lorentzian metric. Then
      any homogeneous manifold $G/\bar H$,  where $\bar H \subset H$
      is a closed subgroup is admissible. We say that $M = G/H$ is
       a { \bf minimal admissible }  manifold (and the  stabilizer $H$ is
       { \bf  maximal admissible} )
      if  there is no   admissible connected  compact Lie subgroup
       $\tilde H$ which contains $H$ properly.\\

   The main goal of the paper is  to  describe  minimal
   admissible manifolds $M = G/H$ of a semisimple Lie group $G$
   and  determine invariant Lorentzian metrics on them.\\

    In section 2, we fix notations and recall  an infinitesimal
    description of  invariant pseudo-Riemannian metrics on a
    homogeneous  manifold $M = G/H$ in terms of the Lie  algebras
    $\gg, \gh$ of the groups $G,H$.\\
    In section 3, we give a necessary and sufficient  conditions
    that a proper  homogeneous manifold  admits an invariant
    Lorentzian metric. We  also  give a description of minimal
     admissible  manifolds $M=G/H$ of a  group $G$ which is a direct product
     $G = G_1 \times G_2$. This reduces  the  classification of
     minimal admissible manifolds of  a semisimple Lie group $G$ to
     the case of a simple  group.\\
     An explicit  description of  minimal  admissible  manifolds $M = G/H$ of
     a simple compact Lie group $G$ and  invariant Lorentzian
     metrics on $M$  is given in section 4. Any  such manifold $M =
     G/H$ is the total space on the canonical $T^1$-bundle
     $$\pi : M = G/H= G/H_{\a} \to F_{\a} = G/H_{\a}\cdot T^1$$
     over  a minimal  adjoint orbit
     $$\Ad_G t_{\a} = G/Z_{G}(t_{\a}) = G/H_{\a}\cdot
     T^1.$$
      The minimal adjoint orbits corresponds  to  simple roots $\a$
      of $G$  and  are the  orbits of  elements $t_{\a}$ of a Cartan
       subalgebra associated with the  corresponding   fundamental
       weights. The  stabilizer $H_{\a}$ is the semisimple part of
       the centralizer $Z_G(t_\a)$. The Dynkin diagram of $H_\a$ is
       obtained  from the Dynkin diagram of $G$  by deleting the
       vertex $\alpha$. Invariant Lorentzian metrics in $M=G/H_{\a}$
       are described in terms of invariant Riemannian metrics in
       $F_\a$  and  the  invariant connection in the  bundle $\pi$. If
       $M$ is not  the  total space of  the sphere bundle over a compact
       rank one  symmetric space, then they
       depends on $m(\a) + 1$ real parameters, where $m(\a) $ is the
       Dynkin mark  associated with the root $\a$.\\

       The section 5 is devoted to investigation  of minimal homogeneous
       Lorentzian manifolds  $M = G/H$ of a  simple noncompact Lie
       group $G$. If $G$ has infinite center, then  the stabilizer
       $H $ is  a  maximal compact   subgroup of $G$.\\
         In the case of  a finite center, the  coset space  $S = G/K$ by
         a maximal compact subgroup $K$   is an  irreducible
         Riemannian symmetric space  with  the  symmetric
         decomposition $\gg= \gk + \gp$.  Let $H\subset K$ be a closed  subgroup
         and
         $$ \gg = \gh + \gm = \gh + (\gn + \gp)$$
          the  corresponding reductive
         decomposition, where  $\gk = \gh + \gn$. The
           subgroup $H$  is admissible if the space $\gm^H = \gn^H + \gp^H$
           of $\Ad_H$-invariant  vectors is nontrivial.  We say that
           the  associated   admissible manifold $M = G/H$ belongs
           to the class I
            if $\gn^H \neq 0$ and  belongs to the class II  if $\gp^H \neq
           0$.\\
            Geometrically, an admissible manifold $M = G/H$ belongs
            to the class I
             if  it admits an invariant Lorentzian metric
             such that the  projection $\pi : M=G/H \to S = G/K$ is
             a pseudo-Riemannian  submersion with  Lorentzian
             totally geodesic fibres $K/H$. In particular, the
             orbits of an invariant time-like  vectors field on $M$
              are circles.
              An admissible manifold $M = G/H$ belongs to the  class II, if
              it admits an invariant Lorentzian metrics with an
              invariant time-like vector field  which  generates a
              noncompact 1-parameter subgroup $\bR$.\\
              Classification of  minimal  admissible manifolds $M =G/H$
              of a simple noncompact Lie group $G$ reduces to
              description of maximal admissible subgroup $H$ of  the
               compact Lie group $K$. This problem had  been solved
               in  section 4.\\
                The classification of  admissible manifolds of
                class II  of a simple Lie group $G$  reduces to
                determination  of the  stabilizers $H= K_v$ of
                minimal orbits  for  the isotropy  representation
                 $$ j : K \to SO(\gp)$$
                 of the symmetric space $S = G/K$. As an example,
                 we determine    such  stabilizers $K_v$  for  the
                 group $SL_n(\bR)$    and  for all simple Lie groups of
                 real rank one and  describe  invariant
                 Lorentzian metrics on  the   associated minimal admissible
                 manifold $M = G/K_v$.\\
                   Starting  from the list of irreducible symmetric spaces
                   $G/K$ of dimension $m \leq 10$, by analyzing
                   the   isotropy representation $j(K)$ we derive
                   also  the list of all class II minimal admissible
                   manifolds $M^d = G/H $ of  dimension $d \leq 11$
                   and  describe  invariant Lorentzian metrics on
                   $M^d$.

\section{Preliminaries}
 By a homogeneous manifold $M=G/H$ we will  understand  the
 homogeneous manifold of a  connected Lie group $G$  modulo a closed
 {\bf connected }  subgroup $H$. We identify  the  tangent  space
 $T_oM$ at the point $o =eH$ with the coset space $V = \gg/\gh$
 where $\gg =\Lie(G)$ is the Lie algebra of $G$  and $\gh = \Lie H$
 is the subalgebra  associated with  the subgroup $H$. We denote  by
 $j : H \to GL(V)$ (resp.,  $j : \gh \to \ggl(V)$)
 the isotropy representation of
the  stability subgroup $H$ ( resp., the stability subalgebra
$\gh$). It is induced by the  adjoint representation of $H$ (resp.,
$\gh$). Since the group $H$ is connected, a tensor $T$ in $V$ is
$j(H)$-invariant if and only if it is $j(\gh)$-invariant, that is
$j(h)T =0$ for all $h \in \gh$.\\
 Recall the  following
 \bp There is a natural bijection between $G$-invariant Riemannian \\
 (resp., Lorentzian) metrics in a homogeneous space $M = G/H$  and
 $j(\gh)$-invariant  Euclidean ( resp., Lorentzian) scalar products
$g_o$  in $V$. An invariant  scalar product $g_o$ defines  the
metric, whose value $g_x$ at a point $x = L_ao :=ao , \, a \in G$ is
given by
$$     g_x := (L_a)^* g_o = g_o((L_a)_*^{-1}\cdot, (L_a)_*^{-1}\cdot).$$
 \ep
Sometimes we will identify $g_o$  and $g$  and say that $g_o$ is  an
invariant metric in $M$.

 Recall that if the group $G$ acts effectively on  a pseudo-Riemannian
   homogeneous
 manifold $M = G/H$, then the isotropy representation is exact  and
 the  stability  subgroup $H$ is isomorphic  to the  isotropy group
 $j(H) \subset GL(V)$. In  particular, we have
\bp A homogeneous manifold $M = G/H$  admits  an invariant
Lorentzian metric  if and only if the isotropy representation  $j$
defines an isomorphism of the stability group $H$  onto a subgroup
$L$  of the connected Lorentz group $SO^0(V)$ or, equivalently,
isomorphism of the stability subalgebra $\gh$ onto a subalgebra
$\gl$ of the Lorentz algebra $so(V)$. \ep

 A homogeneous manifold $M = G/H$ is called to be
 { \bf reductive} if there is  an $\Ad_H$-invariant  (reductive) decomposition
$$  \gg = \gh +\gm .$$
 In this case,  the complementary to $\gh$  subspace $\gm$ is
 identified  with the  tangent space $T_oM = \gg/\gh$  and  the
 isotropy representation is identified with  the  restriction
 $\Ad_H|_{\gm}$ of the  adjoint representation.\\
Any homogeneous manifold with a compact stabilizer is reductive.
\section{Invariant  Lorentzian metrics on a proper  homogeneous $G$-manifolds }

\bd An action of a Lie group $G$ on a manifold $M$ is called
{ \bf  proper} if the map
$$G \times M \to M \times M, \,\, (a,x) \mapsto (ax,x)$$
is proper, or, equivalently, $G$ preserves a complete Riemannian
matric on $M$. In this case $G$-manifold $M$ is called
{\bf proper}.
 \ed
 The orbit  space $M/G$ of a proper $G$-manifold is a metric space  and has
 a structure of a stratified manifold.\\

 For  a nonproper $G$-manifold, the topology of the orbit space can be  very bad,
 for example, non-Hausdorff, see e.g.  the action of the Lorentz group on the
  Minkowski space.\\
 On the other hand, in most cases   the isometry group of a {\bf  compact}
 Lorentzian manifold is compact   and, hence,  acts properly. G. D'Ambra
 \cite{DA} proved that  the isometry group of any  simply connected   compact
  {analitic} Lorentzian manifold  is compact (hence,  it act properly).
  M. Gromov \cite{DAG} states
  the problem of  description of all compact Lorentzian manifolds which admits a
 noncompact (= nonproper) isometry group.  It is a special case of his more general
  problem of  classification of geometric structures of finite order on compact
  manifold  with a noncompact group of automorphisms.
Recall the following
 \bp
 Let $M=G/H$ be a homogeneous manifold with an effective action of $G$.
 Then the following conditions are equivalent:\\
 a)  $M=G/H$ is  proper;\\
 b) the  stabilizer $H$ is compact. \\
c)  $M$  admits  an invariant Riemannian metric (which is  defined
by an $H$-invariant Euclidean metric $g_o$
in $T_oM\, o = eH \in M$)\\
\ep
 An  $H$-invariant metric $g_o$ can be constructed  as the center
of the ball of minimal radius in $S^2(T^*_oM) $ ( w.r.t. some
Euclidean metric $g_1$)  which contains the orbit $j(H)g_1$.

{\color{blue} \large}
\subsection{A criterion for existence of an
invariant Lorentzian metric on a proper  homogeneous manifold $M =
G/H$}

\bp A proper homogeneous manifold $M= G/H$  admits
an invariant Lorentzian metric if and only if 
 the isotropy group $j(H)$ preserves an 1-dimensional subspace $L =
\mathbb{R}v \subset V = \gg/\gh$.\\
 Moreover, let $h$ be a $j(H)$-invariant Euclidean scalar product  and
  $\eta$ is the 1-form  which defines  the hyperplane $L^\perp =
\ker \eta $ orthogonal to  $L$. Then one can associate  with $(L,h)$
an invariant Lorentzian  scalar product
$$   g_0 = h  - \lambda \eta \otimes \eta$$
where $\lambda >0$ is sufficiently big number, which defines an
invariant Lorentzian metric in $M$.
 Any invariant Lorentzian metric
 can be obtained by this construction.
 \ep
\pf The first  claim is obvious.  Now  we prove that  any invariant
Lorentzian metric $g$ on $M$  is obtained by this construction.  The
restriction
  $g_o =g|_o$ is   a $j(H)$-invariant Lorentzian scalar product in $V = T_oM$
and $j(H)$ is a compact subgroup  of the group $SO(V) = SO_{1, n-1}$
which preserves $g_o$. Hence is belongs to a maximal compact
subgroup $O_{n-1} \subset SO_{1,n-1}$ which preserves  a time-like
line  $L= \bR t \in V$. Then
$$h := \lambda \eta \otimes \eta + g_o $$
 for $\eta := g_o(t,\cdot)$ and sufficiently big $\lambda >0$ is a
$j(H)$-invariant Euclidean metric such that
 $g = - \lambda \eta \otimes \eta + h $.
 So the Lorentzian metric $g$ is obtained from a
Riemannian metric (associated with $h$) by the described
construction. \qed

\bc   If $(M=G/H,\,g)$ be a proper  homogeneous Lorentzian manifold
with connected  stabilizer $H$. Then it admits an invariant
time-like vector field $T$ with $g(T,T) =-1$ and the formula
 $$  h = \lambda g\circ T \otimes  g\circ T  + g$$
defines  an invariant Riemannian metric for any $\lambda >1$.
 \ec

We will always assume in the sequel  that the stability
 subgroup $H$ is connected.

\bd A proper homogeneous manifold $M = G/H$  (and   the
corresponding  stability group $H$ )  is called
{\bf admissible} if $M$ admits an invariant Lorentzian metric.\\
  Moreover, a compact subgroup  $H$ is called  { \bf   maximal
  admissible}
  if it is  a maximal compact  subgroup  such that $M=G/H$ admits an
  invariant Lorentzian metric.Then the manifold $M=G/H$ is called a
 {\bf  minimal admissible manifold}.
 \ed

\bc A proper homogeneous manifold $M = G/H$ with a reductive
decomposition  $\gg = \gh + \gm$ is admissible if and only if $\gm^H
\neq 0$ where $\gm^H$ is the space of $\Ad_H$-invariant vectors from
$\gm$. \ec

\bp Any closed subgroup $H'$ of an admissible subgroup $H$  is
admissible. \ep

\pf  Let $\gg = \gh + \gm$ be a reductive decomposition of an
admissible manifold $M = G/H$  and $H'\subset H$ is a subgroup with
$\gh' = {\mathrm Lie} H'$. Then
$$  \gg = \gh' + \gm' = \gh' + (\gp + \gm),$$
where $\gp$ is an $\Ad_{H'}$-invariant  complement to $\gh'$ in
$\gh$, is a reductive decomposition of $G/H'$  and
$$ { \gm'}^{H'} = {\gp}^{H'} + \gm^{H'} \supset \gm^H \neq 0.$$
This  shows  that  $H'$ is  an admissible  subgroup. \qed


 The   above observations  reduce  the  problem of description of  admissible
  homogeneous $G $-manifolds $M = G/H $  to classification of  maximal
  admissible subgroups  $H $ of $G $  and  a description of all  closed subgroup of
  the (compact) maximally  admissible  groups $H$. The problem of construction of
   all invariant Lorentzian metrics on a given  admissible  homogeneous manifold
   $M=G/H $   with a reductive decomposition
 $\gg = \gh + \gm $ reduces to a description of  all invariant Riemannian metrics
  on $M$
 (or , equivalently,  $\ad_{\gh} $-invariant Euclidean scalar products in $\gm $)
  and a description of the  space  $\gm^{H} $  of  $\Ad_{H} $-invariant  vectors
  in $\gm $.

{\bf Example}  Let $M= G/H$  be an admissible  homogeneous manifold
with  a reductive decomposition  $\gg = \gh + \gm$.  Assume  that
the $j(H)$-module $\gm$ admits  a decomposition
$$ \gm = \gm_0 + \gm_1 + \cdots + \gm_k $$
 where $\gm_0$ is a trivial module  and $\gm_i,\, i>0$ are non-equivalent
  irreducible  modules. Then any invariant Lorentzian
 metric on $M$ is defined  by a scalar product of the form
 $$  g = g_0 + \lambda_1 g_1 \cdots + \lambda_k g_k$$
 where $g_0$ is a Lorentzian scalar product, $g_i$ are  invariant  Euclidean
 scalar product in $\gm_i,\, i>0$ and $\lambda_i$ are positive
 numbers.\\

 We will use this construction in the sequel.



\subsection{Minimal homogeneous Lorentzian $G$-manifolds  where $G = G_1 \times G_2$
is a direct product}
 In this subsection we  describe   the structure  of
minimal admissible homogeneous  $G$-manifold  $ M = G/H$ where $G=
G_1\times G_2$ is a direct product of two Lie groups. It reduces the
classification of minimal admissible  homogeneous manifolds of a
semisimple Lie group  $G$  to the case of simple Lie group $G$.\\

The reductive decomposition of  $M =(G_1 \times G_2)/H$ can be
written as
 $$ \gg = \gh + \gm = (\gh_1 + \gh_1 + \gl) + (\gm_1 + \gm_2 + \gl_1 )$$
where  $\gh_i = \gh \cap \gg_i$,  $\gl$ is the complementary to $
\gh_1 + \gh_2$ ideal of $\gh$,  $ \gl_i = \pi_i(\gl) \simeq \gl$ is
the  projection of $\gl$  to $\gh_i$ and  $\gm_i$ is an
$\ad_{\gh}$-invariant complement  to  the compact  subalgebra  $
\gh_i + \gl_i$ in $\gg_i,\, i=1,2$. Assume that the space $\gm_1^H$
of $H$-invariant vectors in $\gm_1$ is not zero.  Then  the
subalgebra $\gh_1 + \gl_1 + \gh_2 + \gl_2$ generates  an admissible
subgroup  which, by maximality of $H$,  coincides with $H$. Hence
$\gl=0$  and  the homogeneous manifold $M$ is  a direct product  $M
= G/H= G_1/H_1 \times G_2 /H_2$. Note that  a subgroup $H_1 \times
H_2 \subset G_1 \times G_2$ is maximal admissible  if one of  the
factors, say $H_1$ is a maximal admissible subgroup of $G_1$  and
the other  factor $H_2$ is  a maximal
compact subgroup of $G_2$.\\

 Assume now that $\gm_i^H =0,\,\, i=1,2$. Then  the  compact
 subalgebra $\gl_1$ must have a  center   and  from  the condition that
 $H$ is a maximal admissible  subgroup  we  conclude  that $\gl_i = \bR
 t_i$  is an 1-dimensional  subalgebra of $\gg_i$ and $\gh_i + \bR
 t_i$ is  its centralizer  in a maximal compact subalgebra  $\gk_i
 $ of $\gg_i$. This  implies the following result.

 \bt \label{Theoaboutproduct} Let $M = G/H$ be a minimal admissible
  homogeneous manifold of a
 Lie group $G = G_1 \times G_2$. \\
 If $H = H_1 \times H_2$ is
 consistent with  the  decomposition of $G$, then one of the
 subgroups $H_1, H_2$, say $H_1$, is  maximal admissible in $G_1$ and the other  subgroup
 $H_2$ is maximal compact subgroup of $G_2$.\\
 If $H$ is not consistent with the decomposition,  then  its Lie
 algebra  has the   form
 $$  \gh = \gh_1 + \gh_2 + \bR(t_1 +t_2)
 $$
 where $\gh_i + \bR t_i = Z_{\gk_i}(t_i)$ is the centralizer of an
  element $t_i \in \gg_i$
 into a maximal compact   subalgebra $\gk_i := \Lie K_i$ of $\gg_i, \, i=1,2$.
 The reductive decomposition  associated with $M=G/H$ can be written
 as
 $$  \gg = \gh + \gm = \gh + (\gm_1 + \gm_2 + \bR(t_1 -t_2))
 $$
where  $\gm_i$ is an $\ad_h$-invariant complement  to
$Z_{\gk_i}(t_i)$ in $\gg_i$.
 \et
 This  theorem can be  applied to the case  when $G$ is a semisimple
  Lie  algebra  and it reduces the description of
 admissible  homogeneous manifolds of a semisimple Lie group $G$  to the
 case of  simple Lie groups.

\section{Homogeneous Lorentzian  manifolds of simple compact Lie group}

  Let $G$ be a compact simple Lie group. The adjoint orbit
   $F = \Ad_G t \simeq G/ Z_G(t)$ of $G$ is called  to be minimal, if   the
    stability  subgroup $Z_G(t)$ ( which is  the centralizer of  an  element
     $t \in \gg$ )  is not contained   properly in  the centralizer of other non-zero
      element $t' \in \gg$. Recall that  the centralizer $Z_G(t)$ is connected.\\
 It is know, see, for example \cite{Al2} that the orbit  $F $  if minimal if and only if $Z_G(t)$ has
 1-dimensional center
 $T^1 = \{ \exp\lambda t \}$  and can be written as $Z_G(t) = H \cdot T^1$ where
  $H$ is a semisimple normal subgroup. Minimal adjoint orbits (up to an isomorphism)
   correspond to
 simple roots $\a$ of the Lie algebra $\gg$. Moreover, the Dynkin diagram of
  the  semisimple  group $H$ is obtained  from the Dynkin diagram of $\gg$ by
   deleting   the vertex  $\a$. We will denote   the minimal orbit associated with
   a simple root $\a$ by $F_{\a}$.
  Below we  give the  list of all such  semisimple subgroups $H$ for
  all simple  Lie groups $G$:\\
 \bean{}
   G &= SU_{n},  \,&  H = SU_p \times SU_q, \, p+q =n,\, p = 1,2, \cdots ,n-1; \\
   G &= SO_n,   \,&  H = SU_p \times SO_q,\,  2p + q =n,\, p = 1,2,\cdots ,[\frac{n}{2}]; \\
   G &= Sp_n,  \,  & H= Sp_p \times Sp_q, \, n=p+q, \, p=1,2, \cdots, n-1;\\
   G &= G_2,  \,  & H= SU_2^{short}, \, SU_2^{long}\\
   G &= F_4,  \, & H= Sp_3, \, SU_3^{short} \times SU_2^{long},\, SU_2^{short} \times SU_3^{long},\, Spin_7 ;\\
   G &= E_6,  \, & H= Spin_{10},\, SU_2 \times SU_5,\, SU_3 \times SU_3 \times SU_2,\, SU_6;\\
   G &= E_7, \, & H = E_6, SU_2 \times Spin_{10}, SU_3 \times SU_5,
    SU_4 \times SU_3 \times SU_2, SU_6 \times SU_2, Spin_{12}, SU_7.\\
   G &= E_8, & H = E_7, SU_2 \times E_6, SU_3 \times
          Spin_{10},SU_4\times SU_5,
     SU_5 \times SU_3 \times SU_2, SU_7 \times SU_2, Spin_{14}.\\
 \eean

 Let $F_{\a} = G/H \cdot T^1$ be a minimal orbit associated with a
 simple root $\a$. Then
 $$\pi: M_{\a} = G/H \to F_{\a}= G/H\cdot T^1 $$
  is a
 principal fibration with the structure group $T^1$. Denote by
 $$\theta : TM_{\a} \to \bR = Lie(T^1)$$
  the $G$-invariant principal connection   defined by  the condition
  $\theta(t) =1,  \,\, \theta(\gp) =0$
where $$  \gg = (\gh + \bR t) + \gp $$ is the reductive
decomposition associated with  the orbit $F_{\a} = G/ H\cdot T^1$.
We say that $\pi$ is the canonical $T^1$ bundle with connection over
the orbit $F_\a$.

It is known that the tangent space $T_oF_{\a}\simeq \gp$  as an
$\Ad_{(H\cdot T^1)}$-module is decomposed  into mutually non
equivalent irreducible submodules

\be \label{decomposition}
 \gp = \gp_1 + \cdots +\gp_m
 \ee

 and the number $m$  of these submodules  equal  to the Dynkin number
  $m({\alpha}) = m_i$ of  the corresponding  simple root $\a = \a_i$  that is the
   coordinate $m_i$ over $\a_i$  in the decomposition   $\mu = \sum_j m_j \a_j $
    of the maximal root $\mu$ with respect to the simple roots
 $\a_1,\cdots,\a_r$. This implies  that any invariant Riemannian metric $g_F$  in
  $F$  at the point
 $o = e(H\cdot T^1)$ is given by
 $$ g_o = \lambda_1 b_1 + \cdots + \lambda_m b_m $$
where $b_j = - B|_{\gp_j} $ is the restriction of the minus Killing form $-B$ to
$\gp_j$  and $\lambda_j$ are positive constants.

\bt \label{TheorMinCompSimple} Any minimal admissible manifold  of a
simple compact Lie group $G$ is the total space
 $M_{\a} = G/H $ of  the canonical fibration  over a minimal orbit
  $F = F_{\a} = G/ H_\a \cdot T^1 $.
Moreover,  if $M= G/H_\a$ is not  the  total space of  the sphere
bundle  of a compact
 rank one symmetric space that is
$$ S(S^n) = SO_{n+1}/SO_{n-1}, Spin_7/SU_3 = S(S^7)= S^7 \times S^6,
S(S^3)= SU_2 \times SU_2/T^1 = S^3 \times S^2;  $$
$$  S(\bC P^n) = SU_{n+1}/ SU_n,\, S(\bH^n)= Sp_{n+1}/Sp_1 \times Sp_{n-2},\, S(\bO P^2)= F_4/Spin_7
$$
then any invariant Lorentz metric $g$ on $M$ is given by
$$   g = -\lambda \theta^2 + \pi^*g_F $$
where $\theta$ is the principal connection, $g_F$ is an invariant
Riemannian metric on $F$ and $\lambda$ is a positive number. In
particular, the metric $g$  depends on $m(\a) +1$ positive
parameters, where $m(\a)$ is the Dynkin mark. \et

\pf Let $M = G/H$ be a minimal  admissible  homogeneous  manifold of
a simple compact Lie group $G$  with  the reductive decomposition
$\gg = \gh + \gm$. Denote by $t \in \gm$  an $\Ad_H$-invariant
non-zero vector. We can assume that $t$ generates a closed
one-parameter subgroup since $H$ preserves pointwise  the  curve
$\exp \lambda t$, hence, also its closure  in $G$. The   centralizer
$\gz(t)$  of $t$ in $\gg$ can be decomposed into a direct sum
$\gz(t) = \tilde \gh + \bR t$ where $\tilde \gh \supset \gh$ is a
subalgebra which generates a closed subgroup $\tilde H$ of $G$.
Since $\Ad_{\tilde H }$ preserves $t$, the homogeneous space
$G/\tilde H$ is admissible  and due to minimality of $M$ it
coincides with $M$.  Hence, $H = \tilde H$  and $Z_G(t) = H\cdot T^1
$  where $T^1$ is the closed subgroup generated by $t$. It is proven
in \cite{AS}, that if $M$ is  not  the total space of the sphere
bundle  of a compact rank one symmetric space,  then  all
irreducible   $(H\cdot T^1)$-submodules  of  the  decomposition
(\ref{decomposition}) remain irreducible and non-equivalent  as
$H$-submodules. This implies the  last claim of the theorem. \qed

\section{ Homogeneous  Lorentzian manifolds of a simple noncompact
Lie group }
 Now  we  consider   minimal admissible  homogeneous
 manifolds  of a  simple noncompact Lie group $G$.
\subsection{Case  when the group $G$ has  infinite center }

 Assume  at first that $G$ has infinite center. It is known that
 such group $G$  acts
 transitively (and  almost effectively)
 on  a non-compact irreducible  Hermitian  symmetric space $S = G/K\cdot \bR $
 with the  symmetric decomposition
  $$ \gg = (\gk + \bR t) + \gp $$
  where  $\bR t$ is  the 1-dimensional  centralizer  of  the Lie  algebra
  $\gk$ of a maximal compact subgroup $K$ of $G$ and
  $\ad_t|_{\gp}$  is $j(K\cdot \bR)$-invariant  complex  structure in the
  tangent space $\gp = T_oS$. Obviously, we get the following
  \bp  Let $G$ be a simple non-compact Lie group and $S = G/K \cdot
  \bR$ the associated  Hermitian symmetric space. Then   the
  manifold $M = G/K$  is the only  minimal admissible  $G$-manifold
  and   all invariant Lorentzian metrics on  $M$  are defined by
  the  scalar product  in $\gm = \bR t + \gp$ of the form
      $$   g= -\lambda \theta^2 + g_{\gp}    $$
  where $\lambda >0$, $\theta$ is the 1-form  dual to  the vector $t$
  (such that   $\theta(t) =0,\, \theta(\gp) =0$) and $\gp$ is  the
  invariant  Euclidean scalar product in $\gp$ which defines  the
  symmetric Riemannian metric in $S$.  In particular,
   $$\pi : M = G/K  \to S = G/K \cdot \bR $$
   is a pseudo-Riemannian submersion.
 \ep
\subsection{Duality }
 Now we  will assume that $G$ is a  simple noncompact Lie group
 with  a finite center. Then  the  quotient $S = G/K$ by a maximal
 compact subgroup $K$ is a symmetric space of noncompact type. We will denote by
  $\hat S = \hat G/K$ the dual  compact symmetric space.
 Let
  $$   \gg = \gh + \gp $$
 be a symmetric decomposition  associated with the  symmetric space  $S$. Then
  the  symmetric decomposition associated with $\hat S$ can be written as
 $$   \hat \gg = \gh + i \gp $$
 where $[iX, iY] = -[X,Y]$  for $X,Y \in \gp$.

 In particular,  the  dual symmetric  spaces  $S, \hat S$  have the  same
 stabilizer $K$   and isomorphic  isotropy representation
 $j(K) = \Ad_K |_{\gp} \simeq \Ad_K |_{i\gp} $. This implies  the  natural
  bijection between  (maximal)
 admissible  subgroups $H \subset K$ of  the dual Lie groups $G$ and $\hat G$.
  In terms of
  homogeneous Lorentzian manifolds  this  can be  reformulated as follows.
 \bp    There exists a natural one-to-one correspondence   between
 proper homogeneous Lorentzian  $G$-manifolds $M = G/H$ of a simple noncompact
  Lie group $G$   and
 homogeneous Lorentzian  manifolds $\hat M = \hat G/H$ of the dual
 compact Lie group $\hat G$
such that  the  stabilizer $H$ belongs to the   subgroup $K \subset
\hat G$. \ep

 \pf Let
 $M= G/H,\, H \subset K$ be an admissible  $G$-manifold  with reductive
 decomposition
 $$  \gg = \gh + \gm := \gh + (\gn + \gp),\,\,   \gk = \gh +\gn$$
  with the invariant Lorentzian metric  defined by an $\Ad_H$-invariant Lorentzian
  scalar product $g_o$ in $\gm = \gh + \gp $,  then  the  dual compact homogeneous
  Lorentzian manifold
  is the  homogeneous  manifold  $\hat M = \hat G/H $  with the reductive
  decomposition

  \be \label{reductivedecomp}
  \hat{\gg} = \gh + \hat{\gm} := \gh + (\gn + i\gp)
  \ee
   and  the metric  defined by the Lorentzian scalar product in $\hat{\gm}$
    which corresponds to  the scalar product $g_o$ under the natural isomorphism
  $$  \hat{\gm} = \gn + i\gp \simeq \gm = \gn + \gp. $$
\qed

\subsection{ A characterization of noncompact homogeneous Lorentzian manifolds
of class I and class II }

  Let $M = G/H, \,\, H \subset K$ be an  admissible homogeneous space of a
   noncompact simple  Lie group $G$ with the reductive decomposition
    (\ref{reductivedecomp}). Then the space
    $\gm^H = \gn^H + \gp^H$ of $j(H)$-invariant vectors is not zero.
 \bd We say  that the admissible  homogeneous  manifold $M= G/H$  belongs to
the  class I if $\gn^H \neq 0$ and belongs to the  class II if
$\gp^H \neq 0$.
 \ed

  Geometrically,  homogeneous spaces of the  class I  and the class II  can be
   characterized  as  follows.
  \bp  An admissible $G$-manifold $M= G/H$ of a simple noncompact Lie group $G$
  belongs to the  class I  if
   it admits an invariant Lorentzian metric such  that $\pi : M = G/H \to S = G/K$
    is a pseudo-Riemannian submersion with totally geodesic Lorentzian
   fibres over  the noncompact  Riemannian  symmetric space $S=
   G/K$. In particular, the  invariant  time-like vector filed
   generate a compact group $S^1$.\\
   An  admissible manifold $M =G/H$ belongs  to the class II if  it  admits an
   invariant Lorentzian metric with a time-like invariant vector
   field, which generates a noncompact 1-parameter  subgroup $\bR$.
  \ep

 \pf Assume that $M$ belongs to the class I.
  Let $t \in \gn^H$ be an $H$-invariant vector  and
 $g = g_{\gn} \oplus g_{\gp}$ an Euclidean scalar product in $\gm$
  which is a sum of $\Ad_{H}$-invariant scalar  product in
 $\gn$ and the unique ( up to a scaling) $\Ad_K$-invariant scalar product
 in $\gp$. Then the invariant Lorentzian metric  in $M$  defined by the
  Lorentzian scalar product of the form
  $ g_{t, \lambda} = g - \lambda g\circ t \otimes g\circ t$
 for sufficiently big $\lambda $ satisfies the  stated property.
\qed

 {\bf Remark } It is possible that  a  minimal  admissible $G$-manifold  belongs
  to the class I   and the class II at the same time. \\

 Let $K \subset GL(V)$ be  a linear Lie group. Recall that  by  the
 (connected)
 stabilizer $K_v$ of a vector  $v \in V$  we understand the {\bf
 connected} component of the  subgroup which preserves   $v$.

 \bd Let $K \subset GL(V)$ be a linear Lie group. The
  orbit $Kv $  of a vector $v \neq 0$ is called  a {\bf minimal orbit} is the
  the ( connected) stabilizer $K_v$ does not  contained properly  in
   the   (connected )  stabilizer  $K_w$
   of any other non-zero vector $w$. Then the stabilizer $K_v$ is called
   a { \bf maximal  stabilizer}.
\ed

 The following obvious proposition  reduces  the   classification of  all minimal
 admissible homogeneous $G$-manifolds $M= G/H$  of  the class I   to  the classification
  of   maximal admissible  subgroups $H$ of the  maximal compact  subgroup $K$ of
  $G$ and the classification of  such manifolds of the class II to  the description
  of  maximal  isotropy subgroups $K_v$ of the isotropy representation
  $\Ad_K|_{\gp}$ of the  symmetric space $S = G/K$.

 \bp  Let $M=G/H$ be a minimal  admissible  homogeneous $G$-manifold of a simple
  noncompact Lie group $G$.
 \begin{enumerate}
    \item[i)] If $M$ belongs to the class I, then  $H$ is  a maximal  admissible
     subgroup  of a maximal compact subgroup $K \supset H$ of $G$.
\item[ii)]  If $M$ belongs to the class II, then $H = K_v$ is a maximal
(connected)
stabilizer of  the isotropy representation  of the Riemannian
symmetric space $S = G/K$.
\end{enumerate}
\ep

Let
  $$\mathfrak{g}= \mathfrak{k} + \mathfrak{p}$$
   be  the symmetric decomposition of   a symmetric space
  $S = G/K$. For any nonzero vector $v \in \mathfrak{p}$ we denote by $\mathfrak{k}_v$ the  stability subalgebra of the isotropy  representation $j(\mathfrak{k})$  and by
  $K_v \subset K$ corresponding connected  stability subgroup.\\
 \bd  The subalgebra $\mathfrak{k}_v \subset \gk$
 (resp., corresponding  subgroup $K_v \subset K $) is called  a
 {\bf maximal stability subalgebra }(resp.,{\bf maximal stability subgroup})
  if it does not contained
  properly in any other  stability subalgebra (resp., stability subgroup) of
  the isotropy representation of $G/K$.
  \ed

 \bp Let $S = G/K$ be a symmetric space of noncompact type   and $H \subset K$
  a maximal  admissible subgroup of $K$ such that  the admissible
   manifold $G/H $ belongs to the class II. Then  $H =H_v$  is a    maximal
    stability subgroup of $K$. Conversely,   any  maximal stability
    subgroup $K_v $ of $K$  is   admissible   and  defines an
    admissible manifold $M=G/K_v$ of the class II.
 \ep

  So the classification of  proper  homogeneous Lorentzian manifolds of a semisimple
 noncompact  group $G$  reduces to description of  maximal stability   subgroups
 $K_v$ of the isotropy representation of the associated  symmetric space $S = G/K$.\\

   Due to theorem \ref{Theoaboutproduct}, it is  sufficient to describe
    such subgroups for  simple Lie groups.\\

\subsection{Homogeneous Lorentzian $ SL_n(\bR) $-manifolds}
In this subsection we classify all minimal  homogeneous Lorentzian
 $G$-manifolds of the class II where  $G= SL_n(\bR)$.\\

  Let $ S = SL_n(\bR)/SO_n $. We identify $S $  with the  codimension one orbit
   $SL_n(\bR)g_0$ of  the Euclidean metric $g_0 \in S^2V^*$ in the
   space $S^2V^*$ of symmetric  bilinear forms
    in $V = \bR^n$ (or  with the space of symmetric matrices).
    In particular,  the  tangent  space
  $T_{g_o}S  = T_o{S}$  is identified  with the space of  $S_0^2(V^*)$
  of traceless (w.r.t. $g_0$) bilinear forms.
   Let $V = U + W$
  be a decomposition of $V$ into a $g_0$-orthogonal  sum of  subspaces of
  dimension $p$ and $q$ ,respectively, and
  $H = SO(U) \times SO(W) = SO_p \times SO_q$ the connected subgroup
  of $SO(V) = SO_n$ which preserves this decomposition. Consider the
  homogeneous manifold
  $$M_{p,q} = G/H := SL_n(\bR)/ SO_p \times SO_q, \,\, p+q =n.$$
   It has
   the natural fibration
$$ M_{p,q} =  SL_n/(SO_p \times SO_q) \to  S = SL_n/SO_n
$$
over the symmetric space    $S = SL_n/SO_n $  with the Grassmannian
  $ Gr_p(\bR^n)= SO_n/SO_p \times SO_q$  as a fibre.
The Grassmannian   is an irreducible  symmetric manifold with the
symmetric decomposition
$$ \gso_n = \gso(V) = (\gso(U) + \gso(W)) + U \wedge W.
$$
Then the reductive decomposition  of the  homogeneous manifold
$$M_{p,q} = SL(V) / SO(U) \times SO(W) = SL_n(\bR)/ SO_p \times SO_q
$$
can be written as

 $$ \gg := \gsl(V) = \gh + \gm = (\gso(U) + \gso(W)) +
  ( \mathbb{R}b + U^* \wedge W^* + S_0^2U^* + S_0^2W^* + U^*\vee W^* )$$
where $\vee$ is the symmetric product,
 $b := q g_0|_U - p g_0|_W$  and $ S_0^2U^*, S_0^2W^*$ are irreducible
submodules of traceless  bilinear forms.
  As a $j(H)$-module, the tangent
space $\gm$  is isomorphic to
$$ \gm = \bR b + (U\otimes V)\otimes \bR^2 + S^2_0 U + S^2_0 W.
$$
In particular,
$$\gm^H =  \bR b   \neq 0.$$

 We get
 \bp The homogeneous manifold  $M_{p,q} $ is an admissible
 manifold. Any invariant Lorentzian metric on it is defined by  the
 scalar product of the form
$$   g= - \lambda_1 b^* \otimes b^* +  g_1\otimes g_{\mathbb{R}^2} + \lambda_2 g_2 +
\lambda_3 g_3
$$
where  $\lambda_i, i=1,2,3 $ are positive constants, $g_1, g_2, g_3$
are the Euclidean scalar products in $ U\otimes  V, S^2_0U $ and $
S^2_0W $ respectively, induced by the metric $g_0$
 and $g_{\mathbb{R}^2}$ is an Euclidean scalar product in $\bR^2$.

 \ep
The  following theorem shows   that the spaces $M_{p,q}$ exhaust all
minimal  homogeneous Lorentzian $SL_n(\bR)$-manifolds of the class
II.

\bt  A minimal admissible homogeneous  $SL_n(\bR )$-manifold $M$ of
class II is isomorphic to  the  manifold
 $ M_{p,q} =  SL_n/(SO_p \times SO_q)$ for some $p,q$ with $p+q=n$.
  \et

\pf The isotropy representation $j$ of the  symmetric space $S=
SL_n(\bR)/SO_n$
 is the standard representation of  $K = SO_n$  in the space
   $T_0S = S_0^2 \bR^n$ of  traceless symmetric matrices.\\
 The  stability subgroups of $j(SO_n)$
are
 $ SO_{p_1} \times \cdots \times SO_{p_s}$
 and
 maximal admissible
subgroups  are  $SO_p \times SO_q$.  They defines  manifolds
$M_{p,q}.$ \qed

\subsection{Homogeneous Lorentzian  $G$-manifolds  where $G$ is a
 simple Lie group of real rank one}
In this  subsection  we describe   minimal  homogeneous Lorentzian
manifolds  $M= G/H$ of the  class II for  all simple  Lie group $G$
of real rank 1. The isotropy group $j(K)$  of the associated  rank
one symmetric space $S= G/K$  acts transitively on the unit sphere
in $T_oS$ and
 the stability subgroups $K_v$  of a point $0 \neq v \in T_oS$ is unique
  ( up to a conjugation), hence, maximal.

The list of all noncompact rank one symmetric  space $S = G/K$ is
given below, see \cite{H}.
\smallskip

{ \bf List of rank one noncompact  symmetric spaces $S= G/K$.}\\
$$\bR H^n= SO^0_{1,n}/SO_n,\,\,\bC H^n = SU_{1,n}/U_n,\,\,\bH H^n =
Sp_{1,n}/Sp_1 \times Sp_n,\,\, \bO P^2 = F_4/\mathrm{Spin}_9.$$
\smallskip

We describe  corresponding minimal admissible  manifolds $M= G/H =
G/K_v$ of the  class II  for each of these  groups together with the
reductive decomposition  $\gg = \gh + \gm$  and   the decomposition
of  the tangent space $\gm$ into irreducible $j(H)$-modules. It
allows to give  an explicit  description of all invariant Lorentzian
metrics on $M$.\\

\subsubsection{Case of the group  $G= SO^0_{1,n}$}
 Let $V = \bR^{1,n}$  is  the Minkowski vector space and
 $V=  \bR e_0 + E$  its
  decomposition  where  $e_0, \, e_0^2 =-1, $ is a unit  time-like vector
   and $E = e_0^\perp $.
The hyperbolic space  is  the orbit
 $\bR H^n = G/K = SO^0_{1,n}e_0$
and $E =T_{e_0}\bR H^n  $ is the tangent space with the  standard
action of the isotropy group $SO_n = SO(E)$. We will identify  the
Lie algebra $\gso_{1,n} = \gso(V)$ with the space $\Lambda^2V$ of
bivectors. Then the  reductive decomposition of $G/K$ is given by
$$
\gg = \gh + \gp = \Lambda^2E + e_0 \wedge E.
$$
The stability subalgebra  $\gh = \gk_{e_1}$ of a unit vector $e_1
\in E$  is  $\gso(W) = \Lambda^2 W $ where $W = e_1^\perp$ is the
orthogonal complement of $e_1$ in $E$. This implies
 \bp  The only class II  minimal admissible
 manifold  of the group $G = SO_{1,n}$ is the manifold $M = SO^0_{1,n}/ SO_{n-1}$.
 It has  the reductive decomposition
 $$ \gso_{1,n} =\gso(V)=  \gso(W) + (\bR( e_0 \wedge e_1) +
 e_0 \wedge W + e_1 \wedge W) $$
where
 $$\bR^{1,n} = V= \bR e_0 +\bR e_1 + W $$
is  an orthogonal decomposition of  the Minkowski space $V$. In
particular, $\gm^H = \bR (e_0 \wedge e_1)$.
  \ep

\subsubsection{Case of the group $ G=  SU_{1,n}$}
 Let $\bC^{1,n} =V$ be  the  complex  pseudo-Hermitian  space  with  the
 Hermitian scalar product $  <.,.>$ of  complex signature $(1,n)$ and
$$  V = \bC e_0 + E  = \bC e_0 + \bC e_1 + W $$
an  orthogonal decomposition, such that
 $$<e_0,e_0> =-1, \, <e_1, e_1> =1.$$
The complex hyperbolic space  is the orbit
 $$\bC H^n = SU_{1,n}[e_0] = SU_{1,n}/U_n
$$
of  the point $[e_0]:= \bR e_0 \in PV $ in the projective space $PV
= \bC P^{n+1}$.  The tangent space  $T_{[e_0]}\bC H^n$ is identified
with  $E = \bC e_1 + W$. In matrix notations (with respect to an
orthonormal basis $e_0, e_1, \cdots, e_n$ of $V$) the reductive
decomposition of $\bC H^n$ can be written as
$$
\gg = \gh + \gp = \gu_n + \bC^n,
$$
$$
\gu_n =
\{
\left(%
\begin{array}{cc}
  -\a & 0 \\
  0 & A \\
\end{array}%
\right)  ,\,\, A \in \gu_n,\,  \a = \tr A \},\,\,\,
 \gp = \{ X:=
\left(%
\begin{array}{cc}
  0 & X^* \\
  X & 0 \\
\end{array}%
\right),\, X \in \bC^n, X^* : = \bar X^t \}.
$$
 The  stability subalgebra  $\gk = \gsu_n \oplus \bR z_0 $, where
 $$  z_0 = i\diag (1, - \frac{1}{n}\Id_{n}). $$
We  identify  the tangent space $\gp = T_{e_0}\bC H^n  = E$ with the
space $\bC^n$ of columns. Then  the  subalgebra $\gsu_n$ acts in
$\gp = \bC^n$ by the matrix multiplication and  $z_0$  as the
multiplication by $-\frac{n}{n-1}i$. \\
The element $v = (1,0, \cdots, 0)^t \in \bC^n = T_{e_0}\bC H^n =\gm$
corresponds to the  matrix
$$ v = e_1 \otimes e_0^* - e_0 \otimes e_1^* =
\left(%
\begin{array}{ccc}
  0 & 1 & 0 \\
  1 & 0 & 0 \\
  0 & 0 & 0_{n-1} \\
\end{array}%
\right).
 $$

 The stabilizer  $H= K_{v} \simeq U_{n-1}$ has  the Lie
algebra
$$ \gh = \gk_{v} = \gsu_{n-1} \oplus \bR z  $$
where $\gsu_{n-1} = \gsu(W)$ acts trivially on $e_0,e_1$  and with
respect to the decomposition $V = \bC e_0 +\bC e_1 + W$  the matrix
$z \in \gh \subset su_{1,n}$ is given by
$$z = i\diag (1,1, -\frac{2}{n-1}\Id_W). $$
The stability  subalgebra $\gh = \gsu_{n-1} \oplus \bR z$
annihilates the 2-dimensional space
$$ \bC v =
\{ cv = ce_1 \otimes e_0^* - e_0 \otimes (ce_1)^*  \}\subset \gm.
$$
The Lie algebra $ \gsu_{n-1} \subset \gh$ acts in the standard way
on the complementary subspace $\gp' = \{ w \otimes e_0^* - e_0
\otimes w^*,\,\, w \in W \} \subset \gp$ isomorphic to $W$. The
element $z$ acts on
$\gp'$ as a multiplication by $-\frac{n+1}{n-1}i$.\\
The reductive decomposition of the sphere  $K/H=  U_n/ U_{n-1}$  has
the form
$$ \gk = \gh + \gn = (\gsu_{n-1} + \bR z)+ (\bR z' + \gn' )
$$
where $z' := \diag(1,-1,0_{n-1})$ and
$$ \gn' := \{ w\otimes e_0^* - e_0 \otimes w^*,\,\, w \in W   \}. $$
The $j(H)$- invariant subspace $\gn^H = \bR z'$  and    $j(z)$ acts
on $\gn' \simeq \bC^{n-1}$ as multiplication by $-\frac{n+1}{n-1}i$.
  We get
   \bp
The only minimal admissible $SU_{1,n}$- manifold of   the class II
is the manifold $M = SU_{1,n}/ U_{n-1}$ with the reductive
decomposition
$$
\begin{array}{cccccc}
 \gsu_{1,n} \,= &  (\gsu_{n-1}\, + \bR z)\, + & (\bR z'\,  +&\gn' + \,& \bC v  &
 +\gp')\\
 j(z)       &   0            &   0    &-\frac{n+1}{n-1}i \,\,\,&0    & -\frac{n+1}{n-1}i
\end{array}
$$
( We indicate the action of the central  element $z \in \gh$ on the
corresponding irreducible subspaces.) \ep

Since $\gn^H = \bR z' \neq 0$,   the   manifold  $M$ belongs also to
the the class I.
 The next proposition, which follows from  Theorem
 \ref{TheorMinCompSimple} and
 Theorem \ref{Theoaboutproduct},
   describe  all minimal admissible  $SU_{1,n}$-manifolds of
 the class I.
Let $gu_n = \bR z_0 + \gsu_n$ be the Lie algebra of the group $U_n$
and $ a \in \gsu_n $ an element  such that $\bR(z_0 + a) $ generate
a closed subgroup $T^1_a$ of $U_n$.

 \bp Any class I minimal admissible  $SU_{1,n}$-manifold  is isomorphic to
 one of the manifolds :\\
 $ a)\,\,  SU_{1,n}/ SU_n,\,\\
 b)\,\, SU_{1,n}/ T^1_a \cdot Z_{SU_n}(a),
 \,\,\,  0 \neq a \in \gsu_n $ or \\
$c)\,\,  SU_{1,n}/ T^1_0 \cdot H' $ where $H'$ is
 a maximal admissible subgroup of $SU_n$.
\ep
 \pf We have to describe maximal admissible subgroups $H$ of $ U_n$.
If the Lie algebra $\gh$ of  $H$  contains  the  center $\gz = \bR
z_0$,    we get c). If the projection of $\gh$ on $\gz$ is trivial,
then $\gh = \gsu_n$  and we get  a). If the projection is non
trivial, then  $\gh = \bR(z+a) \oplus \gh'$ for some non-zero $a \in
\gsu_n$, where $\gh' $ is a subalgebra of $\gsu_n$. The reductive
decomposition  of $\gu_n$  can be written as
$$    \gu_n = \gh + (\bR z + \gm')$$
 where $ \gsu_n = (\bR a + \gh')+ \gm'$ is a reductive decomposition of
 $\gsu_n$. The  maximally of $\gh$ implies that
  $\bR a + \gh' = \gz_{\gsu_n}(a)$ and we  get b),  where $T^1_a$ is
  the 1-parameter  subgroup generated by $z + a$.
\qed

\subsubsection{Case of the group  $G = Sp_{1,n}$}
Let $V = \bH^{1,n}$ be  the quaternionic vector space with a
Hermitian form  $<.,.>$ of quaternionic signature $(1,n)$ and
$$  V  = \bH e_0 + E = \bH e_0 + \bH e_1 + W$$
its orthogonal decomposition  with $<e_0,e_0> = - <e_1, e_1> =-1$.
 The quaternionic  hyperbolic space $\bH P^n = G/K =  SU_{1,n}/ Sp_1 \cdot
Sp_n$ is the orbit $ \bH H^n = SU_{1,n}[e_0]$ in the  quaternionic
projective space $\bH P^{n+1}$. The tangent space $T_{[e_0]}\bH H^n
=E$. In terms of  an orthonormal basis $e_0, e_1, \cdots , e_n$ of
$\bH^{1,n}$, the  reductive decomposition of $ \bH H^n  $ is given
by
$$  \gsp_{1,n} = \gh  + \gp
$$
$$
 \gh =
 \{ \left(
\begin{array}{cc}
  \a & 0 \\
  0 & A \\
\end{array}\right) ,\, \a \in {\mathrm Im} \bH = \gsp_1,\, A \in \gsp_n
 \},\,\, \gp = \{
\left(%
\begin{array}{cc}
  0 & X^t \\
  X & 0 \\
\end{array}%
\right), \,\, X \in \bH^n \}.
$$
 Under identification $T_{[e_0]}\bH H^n = E = \gp $,  the vector
 $e_1$  is identified with  the matrix
 $$ v = e_1 \otimes e_0^* - e_0 \otimes e_1^* =
 \left(%
\begin{array}{ccc}
  0 & 1 & 0 \\
  1 & 0 & 0 \\
  0 & 0 & 0_{n-1}
\end{array}%
\right).
$$

 The stabilizer $H= K_{v}  $  of  the vector $ v = e_1 \in E=
T_{[e_0]}\bH H^n $  has the Lie algebra
$$
\gh =\gsp_1 + \gsp_{n-1} =  \{ (\a, A):=
\left(%
\begin{array}{ccc}
  \a & 0 & 0 \\
  0 & \a & 0 \\
  0 &  0 & A \\
\end{array}%
\right)  ,\, \a \in {\mathrm Im} \bH, \,  A \in \gsp_{n-1} \}.
$$
The action $j(\a, A)$  on the space
$$  \gp = \bH v + \gp' =  \{
(x,X): =
\left(%
\begin{array}{ccc}
 0 & x^* & X^* \\
  x & 0 & 0 \\
  X & 0 & 0 \\
\end{array}%
 \right)
,\,  x \in \bH, X \in \bH^{n-1}, X^* = \bar X^t   \}
$$
 is given by
$$  j(\a,A)(x,X) = (\a x - x \a,  AX - X\a).
$$
 The  complementary subspace $\gn$ to $\gh$ in $\gk$ is given  by
 $$
 \gn = {\mathrm Im}H + \gn' =
   \{
(y',Y): =\left(
\begin{array}{ccc}
  -y' & 0 & 0 \\
  0 & y' & -Y^* \\
  0 & Y & 0_{n-1} \\
\end{array}%
\right) ,\, y' \in {\mathrm Im}\bH, \, Y \in \bH^{n-1} \}.
 $$
The action  $j(\a,A) \in j(\gh)$ on $(y',Y) \in \gn $ is given by
$$
j(\a,A)(y',Y) = (\a y' - y' \a, AX - X \a).
$$

These formulas  implies   the following proposition.

 \bp
The  minimal admissible $SU_{1,n}$-manifold of the  class II is the
manifold  $M = Sp_{1,n+1}/ Sp_1 \times Sp_{n-1} $ with  the
reductive decomposition
$$
\gsp_{1,n} = \gh + \gm = (\gsp_1 + \gsp_{n-1})+ ({\mathrm Im }\bH +
\gn' +  \bR v  + ({\mathrm Im}\bH) v  + \gp'  ).
$$
where $\gn' \simeq \gp' \simeq \bH^{n-1}$.\\
 In particular,  the
space $\gm^H = \bR v$ is one-dimensional. As $(\gsp_1 +
\gsp_{n-1})$-module, the tangent space $\gm$ is isomorphic to
$$ \gm = \bR v + \gsp_1 \otimes \bR^2 + \bH^{n-1}\otimes \bR^2
$$
with the natural action of $\gh= \gsp_1 + \gsp_{n-1}$. Any invariant
Lorentzian metric in $M$ is defined by   the scalar product of the
form
$$g = -\lambda v^* \otimes v^* + g_1 \otimes h_1 + g_2 \otimes
h_2$$
 where $g_1, g_2$ are   invariant  Euclidean scalar products  on $\gsp_1,
 \bH^{n-1}$,  respectively,  $h_1,h_2$ are any invariant Euclidean
 scalar products in $\bR^2$,  and $\lambda$ is a positive constant.
 \ep
\subsubsection{Case of the group  $G = F_4$}
We consider  the noncompact exceptional Lie group $F_4$ with maximal
compact subgroup $K = Spin_9$. The  symmetric space $\bO H^2 =G/K=
F_4/Spin_9$ is    dual to the octonion plane. The isotropy group
$j(K)$ acts transitively on the unit sphere $S^{15}$ in the tangent
space $T_0 \bO H^2 = \gm $ with stability subgroup $Spin_7$.  The
irreducible spinor  $Spin_9$-module  $\gp \simeq \bR^{16}$ as a
$Spin_7$-module
  is decomposed into  the following  irreducible  $Spin_7$-submodules
   $$\gm =  \bR v  + \gm_1^8 + \gm_2^7 $$
   where  $\gspin_7 + \gm_2^7 \simeq \gspin_8 \simeq \gso_8$ and   $\gm_1 $
   is $8$-dimensional  spinor $Spin_7$-module.
     We get
   \bp  The minimal admissible $F_4$-manifold is the manifold $ M = F_4/ Spin_7
   $  with the  reductive decomposition
   $$ \gf_4 = \gspin_7 +\gm = \gspin_7 + (\bR v +\gm_1^8 + \gm_2^7 ). $$
   Any invariant Lorentzian metric is  given by
   $$  g = -\lambda_0 v^* \otimes v^* + \lambda_1 g_1 + \lambda_2g_2$$
   where  $g_1, g_2$ are some fixed Euclidean invariant scalar products in
   $\gm_1^8$ and $\gm_2^7$ and $\lambda_i >0,\, i=0,1,2,$.
   \ep

\section{ Homogeneous Lorentzian  class II manifolds of dimension $d \leq 11$
of a simple noncompact Lie group}

 Here  we describe noncompact  minimal admissible class II  manifolds
 $M = G/H$
of dimension $d \leq 11$ with a simple Lie group $G$. The stability
subgroup $H$  is  the  stability  subgroup  $H= K_v$ of a minimal
orbit  $j(K)v$   of  the  isotropy representation
$$j : K \to GL(\gp)$$
 of the  corresponding  noncompact symmetric space $S = G/K$ of
 dimension  $m \leq 10$.
Since  we  already treated  the case of $G= SL_n(\bR)$ and the case
of  real rank one, it is  sufficient to consider simple Lie groups
$G \neq
SL_n(\bR)$ of real rank greater then one.
Any such manifold $G/H$  admits a fibration over  a noncompact symmetric space of
dimension $m \leq 10$. \\
\smallskip
Due to section 5.4, we may assume that $G \neq SL_n(\bR)$.

 { \bf List
of symmetric spaces $S = G/K$ of dimension $m \leq 10$, where $G
\neq SL_n(\bR)$ is a simple noncompact group of real
rank $>1$   ( up to a local isomorphism )}\\
\smallskip

$$
\begin{array}{llc}
   AIII, m=8  &  Gr^8_2(\bC^4) = SU_{2,2}/S(U_2 \times SU_2), & \gp =
\bC^2 \times \bC^2 \\
  BDI,p=2, q=3,4,5 & Gr^{2q}_2(\bR^{2+q})= SO_{2,q}/SO_2 \times SO_q & \gp= \bR^2 \otimes \bR^q \\
  BDI,p=3, q=3 &  Gr^{9}_3(\bR^{6})= SO_{3,3}/SO_3 \times SO_3                             &  \gp = \bR^3 \otimes \bR^3 \\
G, m=8    & G_2/SU_2 \times SU_2 & \gp = \bC^2 \otimes \bC^2.
\end{array}
 $$
{\bf Remark} Here  we take into account  the local  isomorphism  of
the  following   symmetric spaces : \\
$ SU_{1,1}/U_1 \simeq
SO_4^*/U_2
\simeq Sp_1(\bR)/U_1 \simeq SL_2(\bR)/SO_2 = \bR H^2, $\\
$  Sp_{1,1}/Sp_1 \times Sp_1 \simeq SO_{1,4}/SO_4 = \bR H^4, $\\
$SO^*_6/U_3 \simeq SU_{1,3}/U_3 = \bC H^3,$\\
$  Sp_2(\bR)/U_2 \simeq  SO_{2,3}/SO_2 \times SO_3 .$\\

Recall  that   local isomorphism means  the isomorphism of  the
universal covering and we   consider all homogeneous   spaces up to
a covering.

\subsection{ Case of the group  $G = SO_{p,q} $}
The isotropy representation of the symmetric space $SO_{p,q}/ SO_p
\times SO_q$ is the  standard  representation of $K = SO_p \times
SO_q = SO(U) \times SO(W), \,\, U = \bR^p,\, W = \bR^q $ in the
space $V = \gp = U \otimes W $.  Any    element $v \in V$  belongs
to the $K_v$-invariant subspace $U(v) \times W(v)$ where
$$ U(v):= i_{W^*}v,\,\,  W(v) = i_{V^*}v $$
are supports  of $v$. Note that $\dim U(v) = \dim W(v) =r $, where
$r$ is the rank of $v$. This reduces  the classification of
$K$-orbits in $V$  to the case when $\dim U = \dim V =r$, that is to
the classification  of  the orbits of  nondegenerate $r \times r$
matrices $v \in {\mathrm Mat}_r$  with respect  to  the natural
action of the group $K= SO_r \times SO_r$.
 Since  any matrix can be decomposed into a product of an orthogonal
 matrix and  a symmetric matrix  and   any symmetric matrix is
 conjugated  by  element from $SO_r$ to a diagonal matrix,  we get
\bl  Any   $K = SO_r \times SO_r$-orbit in the space ${\mathrm
Mat}_r$ contains  a diagonal matrix.  The  orbit of a nondegenerate
matrix is minimal  if  it is  the orbit  of
 the  diagonal matrix  of the form $\lambda D_k$, where
  $$ D_k = \diag (\Id_{r-k}, -\Id_k).$$
The  stability  subgroup   of  the identity matrix $D_0$ is the
diagonal  subgroup  $K_{D_{0}} = SO_r^{diag} \subset K = SO_r \times
SO_r$.  The   stability subgroup $K_{D_{k}} \simeq SO_r$ is a
twisted   diagonal subgroup of $K$  with the Lie algebra
$$
\gk_{D_{k}} =\{ \left(
\left(%
\begin{array}{cc}
  A_{11} & A_{12} \\
  -A_{12}^t & A_{22} \\
\end{array}%
\right),\,
\left(%
\begin{array}{cc}
  -A_{11} & A_{12} \\
  -A_{12}^t & -A_{22} \\
\end{array}%
\right),\right) \}
$$
\el

Using this lemma, one can easily  describe all class II  minimal
admissible manifolds $M^m = SO_{p,q}/H$ of dimension $m \leq 11$. To
state the final result, we  fix some notations.

We  denote by $e_i, \, i=1, \cdots ,p$ an orthonormal basis of $U =
\bR^p$ and by $f_1 ,\cdots, f_q$  an orthonormal basis of $W =
\bR^q$ and we use the identifications
 $$\gso_p = \gso(U)= \Lambda^2U,\,\, \gso_q = \gso(W)= \Lambda^2W.
$$
 Now we  describe  the  minimal admissible  manifolds $M= SO_{p,q}/H = SO_{p,q}/K_v$  associated
 with    minimal orbits $j(K)v$ of different diagonal  elements $v \in V = U\otimes
 W$. We indicate also  the  stability subalgebra
 $\gh= \gk_v \subset \gso(U) + \gso(W)$
 and the  reductive  decomposition
 $$   \gso_{p,q} = \gh + \gm = \gh + (\gn + \gp ) $$
and the subspace $\gm^H$ of invariant vectors. We set
$$ U' = e_1^\perp ,\,\, \, W' = f_1^\perp,\,\,\,  U'' = \span(e_1,e_2),\,
W'' = \span(f_1,f_2),$$
 $$ E = \span (e_1,e_2),\,  F = \span(f_1,f_2)
.$$
a) $ v = e_1\otimes f_1$.\\
     $H= K_v = SO(U')\times SO(W'), $
$$
\begin{array}{lll}
  \gh & = & \gso(U') + \gso(W') \\
  \gn &  =  & (e_1\wedge U' + f_1 \wedge W') \\
   \gp & = & (\bR v  + e_1 \otimes W' + U' \otimes f_1 + U'\otimes W' ) \\
  \gm^H & = & \gp^H = \bR v.  \\
\end{array}
$$

b) $ v = e_1\otimes f_1 \pm e_2 \otimes f_2$,\\
   $\quad  K_v = SO_2^{diag}\times SO(U'')\times SO(W'').$\\
$$
\begin{array}{lll}
     \gh &=  & \bR(e_1 \wedge e_2 \pm f_1 \wedge f_2)+ \gso(U'') + \gso(W''), \\
     \gn &= & \bR(e_1 \wedge e_2 \mp f_1 \wedge f_2)+ E \wedge W'' + U'' \wedge F,  \\
     \gp &= & \bR v  + \bR(e_1\otimes f_2 \mp e_2 \otimes f_1)+ \span(e_1
     \otimes f_2 \pm e_2 \otimes f_1, e_1\otimes f_1 \mp e_2 \otimes f_2)+\\
          &  &  E\otimes W'' + U'' \otimes F + U'' \otimes W''\\
      \gn^H &=  & \bR(e_1 \wedge e_2 \mp f_1 \wedge f_2)\\
      \gp^H &= & \bR v
\end{array}
 $$

c)  $v_{\pm} = e_1 \otimes f_1 +
  e_2 \otimes f_2 \pm e_3 \otimes f_3$.\\
   We assume  for simplicity that  $p=q=3$.\\
 $ K_{v_{\pm}}= SO_3^{diag}\subset K = SO_3 \times SO_3,$\\
 $$
\begin{array}{lll}
\gh &= & \gk_{v_+} = \span(e_i \wedge e_j + f_i\wedge f_j , \, i,j =1,2,3,) \\
\gn &= & \span(e_i \wedge e_j - f_i\wedge f_j ),\\
\gp &= & \bR v_+ + \gsl_3(\bR) = \bR v_+ + \Lambda^2(\bR^3) + S^2_0(\bR^3). \\
\gm^H & = &  \gp^H = \bR v.\\
\end{array}
$$

{\bf Remark}  i) The group $K_{v_+} = SO_3$ acts in the space $\gp =
{\mathrm Mat}_3 = \ggl_3(\bR)$ by conjugation and its preserves  the
1-dimensional space  $\bR v_+$ of scalar matrices  and  acts
irreducibly on the   space $\Lambda^2(\bR^3)$ of skew-symmetric
matrices and on the space $S^2_0(\bR^3) $ of traceless symmetric matrices.\\
ii) The case of the minimal orbit of  the vector $v_-$ is  similar,
but  the description of the  reductive decomposition is more
complicated and it is omitted.

\bp All class II minimal  admissible $SO_{p,q}$-manifolds  $M =
SO_{p,q}/K_v$ of dimension $m \leq 11$  belong to the following
list:
$$
\begin{array}{ccc}
  M^5 = & SO_{2,2}/SO_2^{diag} & v = e_1 \otimes f_1 + e_2 \otimes f_2\\
  M_1^5 =  &  SO_{2,2}/\{e \}\times SO_2 & v = e_1 \otimes f_1 - e_2 \otimes f_2 \\
M^9 =   & SO_{2,3}/\{e\} \times SO_2 , \,& v = e_1 \otimes f_1 \\
 M_1^9 =  &  SO_{2,3}/ SO_2^{diag} & v = e_1\otimes f_1 \pm e_2 \otimes f_2.  \\
\end{array}
$$
\ep
 \pf The proof follows from    given above  description  of   the stability
subgroup $K_v$ of diagonal  elements of the  form
  $$ v = e_1 \otimes f_1, e_1 \otimes f_1 \pm e_2 \otimes f_2,
   e_1 \otimes f_1 + e_2 \otimes f_2 \pm e_3 \otimes f_3 $$
and calculation of the dimension of the corresponding manifold
$SU_{p,q}/K_v$. \qed
\subsection{Case of the group  $G= G_2$}
 The isotropy  action of the  symmetric space $G_2/SU_2 \times SU_2$
 is the standard action of $K = SU_2 \times SU_2$ in the space
 $\gp = \bC^2 \otimes \bC^2 = \ggl_2(\bC)$ of complex matrices.
  The  manifold $ M =
 G_2/K_v$ has dimension $\leq 11$ if $\dim K_v \geq 3$. There is the
 only one such stability subgroup, the  diagonal subgroup
 $SU_2^{\diag}$, which is the  stabilizer of the identity matrix. The
 group $ SU_2^{\diag} $  acts irreducibly on the  subspace
  $ {\mathrm Herm}^0_2 \subset \ggl(\bC)$
 of Hermitian matrices with zero trace   and on the space
 $ i {\mathrm Herm}^0_2(\bC) = \gsu_2$ of skew-Hermitian matrices. We get
 \bp  The only class II  minimal admissible $G_2$-manifold is the
 manifold $ M^{11} = G_2/ SU_2^{diag} $.  It has the following
 reductive decomposition
$$   \gg_2 = gsu_2^{\diag} + (\gsu_2^{adiag}  + \bC \Id  +{\mathrm Herm}_2^0 +
 i{\mathrm Herm}_2^0)
$$
where $\gsu_2^{adiag} $ is  the anti-diagonal subspace, such that
$$  \gsu_2 + \gsu_2 = \gsu_2^{diag} + \gsu_2^{adiag}. $$
 In particular, $\gm^H = \bC \Id \simeq \bR^2$ and
$\gsu_2^{diag}$-module $\gm \simeq \bR^2 + 3 \gsu_2$.
 \ep
\subsection{The main  theorem}
Combining  all obtained results, we get the following  theorem.

\bt All  minimal admissible  class II  manifolds $M^d=G/H$ of
dimension $d \leq 11$  where $G$ is  a simple noncompact Lie group
are described in the Table I. There are also indicated  the maximal
compact subgroup $K$ of $G$ and the  space $\gm = T_oG/K$ of its
isotropy representation, the dimension $m$ of the symmetric space
$G/K$ and the fibre $K/H$ of the natural $G$-equivariant fibration $
M = G/H \to S = G/K$ over the symmetric space $S = G/K$. \et

\centerline {\bf Table I.}
$$
\begin{array}{|c|c|c|c|c|c|}
  \hline
   &  &  &  &  &  \\
 d & M^d & K & \gm & m & K/H \\
 \hline
  \hline
  3 & SL_2(\bR) & SO_2 & \bR^2 & 2 & S^1  \\
 \hline
  5 & SO_{1,3}/SO_2 & SO_3 & \bR^3 & 3 & S^2 \\
 \hline
  7 & SL_3(\bR)/SO_2 & SO_3 & S^2_0(\bR^3) & 5 & S^2 \\
 \hline
  7 & SU_{1,2}/U_1 & U_2  & \bC^2 & 4 & S^3 \\
 \hline
  7 & SO_{1,4}/SO_3 & SO_4 & \bR^4 & 4 & S^3 \\
 \hline
  9 & SO_{1,5}/SO_4 & SO_5 & \bR^5 & 5 &  S^4 \\
 \hline
  11 & SU_{1,3}/U_2 & U_3 & \bC^3 & 6 & S^5 \\
 \hline
  11 & SO_{1,6}/SO_5 & SO_6 & \bR^6 &  6 & S^5 \\
 \hline
  11 & G_2/SU_2^{\diag} & SU_2 \times SU_2 & \bC^2 \otimes \bC^2 & 8 & S^3 \\
  \hline
   \hline
\end{array}
$$

\end{document}